\documentclass [a4paper,twoside,12pt]{article}
\usepackage{amsmath,amsfonts,amssymb}
\usepackage{vmargin,graphicx,theorem}
\usepackage[english]{babel}
\usepackage{enumerate}
\usepackage{color}

\setpapersize[portrait]{A4}
\setmarginsrb{2.5cm}{2.5cm}{2.5cm}{2cm}{0cm}{0.1cm}{0.5cm}{2cm}
\usepackage[all]{xy}

\newcommand\Ex[1][]
{\noindent\textbf{Exemple}{\ifx\relax#1\else\ \emph{(#1)}\fi}. \\ }
\newcommand\Exs[1][]
{\noindent\textbf{Exemples}{\ifx\relax#1\else\ \emph{(#1)}\fi}. \\ }
\newcommand\Exo[1][]
{\noindent\textbf{Exercice}{\ifx\relax#1\else\ #1\fi}. \\ }
\newcommand\Rem[1][]
{\noindent\textbf{Remarque}{\ifx\relax#1\else\ \emph{(#1)}\fi}. \\ }
\newcommand\Rems[1][]
{\noindent\textbf{Remarques}{\ifx\relax#1\else\ \emph{(#1)}\fi}. \\ }

\newtheorem{thm}{Th\'eor\`eme}[section]

\newtheorem{ex}[thm]{Example}

\newtheorem{cor}[thm]{Corollaire}

\newcommand \beq {\begin{equation}}
\newcommand \eeq {\end{equation}}
\newcommand{\R}{{\mathbb {R}}}%reels
\newcommand{\N}{{\mathbb {N}}}% entiers
\def\E{{\mathbb{E}}}%esperance
\def\cE{{\cal E }}

\def\cL{{\cal L }}

\def\I{{\bf  1}}
\def\bpm{\begin{pmatrix}}
\def\epm{\end{pmatrix}}

\newtheorem{theoreme}{Theorem}[section]
\newtheorem{lemme}[theoreme]{Lemma}
\newtheorem{definition}[theoreme]{Definition}
\newtheorem{proposition}[theoreme]{Proposition}
\newtheorem{corollaire}[theoreme]{Corollary}

\newenvironment{exemp}{\noindent{\bf Example. --- }}{\par}
\newenvironment{exemps}{\noindent{\bf Examples}\benum}{\eenum\par}
\newtheorem{rmq}[theoreme]{Remark}
\newtheorem{rmqs}[theoreme]{Remarks}

\newenvironment{preuve}{\noindent{\it Proof. --- }}
{\hfill\rule{1.3mm}{2mm}\par} 
\newenvironment{application}{\noindent{\bf Application. --- }}{\par}
\newenvironment{applications}{\noindent{\bf Applications. --- 
}\benum}{\eenum\par}

\newcommand{\bex}{\begin{exemp}}
\newcommand{\eex}{\end{exemp}}
\newcommand{\bexs}{\begin{exemps}}
\newcommand{\eexs}{\end{exemps}}

\newcommand{\bprop}{\begin{proposition}}
\newcommand{\eprop}{\end{proposition}}

\newcommand{\bthm}{\begin{theoreme}}
\newcommand{\ethm}{\end{theoreme}}

\newcommand{\bcor}{\begin{corollaire}}
\newcommand{\ecor}{\end{corollaire}}

\newcommand{\blem}{\begin{lemme}}
\newcommand{\elem}{\end{lemme}}

\newcommand{\beqna}{\begin{eqnarray}}
\newcommand{\eeqna}{\end{eqnarray}}
\newcommand{\beqnas}{\begin{eqnarray*}}
\newcommand{\eeqnas}{\end{eqnarray*}}

\newcommand{\bmat}{\begin{pmatrix}}
\newcommand{\emat}{\end{pmatrix}}

\newcommand{\bpf}{\begin{preuve}}
\newcommand{\epf}{\end{preuve}}

\newcommand{\benum}{\begin{enumerate}}
\newcommand{\eenum}{\end{enumerate}}

\newcommand{\bitem}{\begin{itemize}}
\newcommand{\eitem}{\end{itemize}}

\newcommand{\brmq}{\begin{rmq}}
\newcommand{\ermq}{\end{rmq}}

\newcommand{\brmqs}{\begin{rmqs}}
\newcommand{\ermqs}{\end{rmqs}}

\newcommand{\bapp}{\begin{application}}
\newcommand{\eapp}{\end{application}}

\newcommand{\bapps}{\begin{applications}}
\newcommand{\eapps}{\end{applications}}

\newcommand{\bdefi}{\begin{definition}}
\newcommand{\edefi}{\end{definition}}

\newcommand{\bequa}{\begin{equation}}
\newcommand{\eequa}{\end{equation}}

\newcommand{\dR}{\ensuremath{\mathbb{R}}}

\theoremstyle{definition}

\newcommand{\ABS}[1]{\ensuremath{{\left| #1 \right|}}} % |1|
\newcommand{\BRA}[1]{\ensuremath{{\left\{#1\right\}}}} % {1}
\newcommand{\NRM}[1]{\ensuremath{{\left\Vert #1\right\Vert}}} % \|1\|
\newcommand{\SBRA}[1]{\ensuremath{{\left[#1\right]}}} % [1]

\newcommand{\PAR}[1]{\ensuremath{{\left(#1\right)}}} % (1)

\begin{document}

\title{Weighted Nash Inequalities}

\author{Dominique Bakry\thanks{Institut de Math\'ematiques de Toulouse, UMR CNRS 5219, Universit\'e Paul-Sabatier}~\thanks{Institut Universitaire de France.} 
 ,  Fran{\c c}ois Bolley\thanks{Ceremade, UMR CNRS 7534, Universit\'e Paris-Dauphine.}, 
  Ivan Gentil\footnotemark[3] \,and Patrick Maheux\thanks{Universit\'e d'Orl\'eans}}

\maketitle

\begin{abstract}
Nash or Sobolev inequalities are known to be equivalent to ultracontractive properties of Markov semigroups, hence to uniform bounds on their kernel densities. In this work we present a simple and extremely general method, based on weighted Nash inequalities, to obtain non-uniform bounds on the kernel densities. Such bounds imply a control on the trace or the Hilbert-Schmidt norm of the heat kernels. We illustrate the method on the heat kernel on $\dR$ naturally  associated with the  measure  with density $C_a\exp(-|x|^a)$, with $1<a<2$, for which uniform bounds are known not to hold.
\end{abstract}

\bigskip


\bigskip


\section*{Introduction}
\label{intro} The classical Nash inequality in $\R^n$ may be stated as 
\beq 
\label{Nash}
 \|f\|_2^{1+n/2}\leq C_n\|f\|_1\| \nabla f\|_2^{n/2}
 \eeq
for all smooth functions $f$ (with compact support for instance) where the norms are computed with respect to the Lebesgue measure. This inequality has been introduced by J. Nash in 1958 (see \cite{nash}) to obtain regularity properties  on the solutions to parabolic partial differential equations. The optimal constant $C_n$ has been computed more recently in~\cite{carlen-loss}. 

 In the more general setting of a symmetric Markov semigroup $(P_t)_{t\geq0}$  one has to replace $\| \nabla f\|_2^2$ by the Dirichlet form $\cE(f,f)$ associated with its generator.  Inequality~\eqref{Nash} implies smoothing  properties of the Markov semigroup in the following way : given a function  $f$, then $\varphi(t)=\| P_tf\|_2^2$  has derivative $\varphi'(t)=-2\,\cE(P_t f,P_t f)$,  so,  by the Nash inequality~\eqref{Nash},
$$
\varphi(t)^{1+n/2}\leq C_n^2\| P_t f\|_1^2\PAR{-{\varphi'(t)}/{2}}^{n/2}  \leq C_n^2\|f\|_1^2\PAR{-{\varphi'(t)}/{2}}^{n/2}.
$$
 Integrating leads to the first bound $\| P_tf\|_2\leq C' t^{-n/4}\|f\|_1$ for $t>0$ and  then to  $\| P_tf\|_\infty\leq C' t^{-n/4} \|f\|_2$ by duality and symmetry of the semigroup. This finally   implies  the classical uniform bound
\begin{equation}
\label{eq-ultra0}
\| P_tf\|_\infty\leq {C'}^2 \, t^{-n/2} \, \|f\|_1
\end{equation}
for $t>0$ by semigroup properties. In turn this implies uniform bounds  on the kernel density of the semigroup such as 
\begin{equation}\label{ptunif}
\vert p_t(x,y) \vert \leq {C'}^2 t^{-n/2}
\end{equation}
for all $x,y$ and $t>0$.

\bigskip

Depending on whether the reference measure is finite or not, Nash inequalities take the general form 
\beq 
\label{Nash2} \|f\|_2^{1+n/2}\leq \|f\|_1[a \, \cE(f,f) + b \, \|f\|_2^2]^{n/4},
\eeq  
where $n$ no longer needs to be an integer. They are one of the many forms of the celebrated Sobolev inequality  
\beq 
\label{Nash3} \|f\|_{2n/(n-2)}\leq a \, \cE(f,f) + b \, \|f\|_2^2
\eeq 
for $n>2$, see \cite{BCLS,saloff11}.  Up to constants, these inequalities are all equivalent to the ultracontractive bound
\beq 
\label{ineq:ultrac}
\| P_tf \|_\infty \leq Ct^{-n/2}\|f\|_1, \quad 0 < t \leq 1
\eeq 
on the Markov semigroup associated to the Dirichlet form $\cE$, hence to {\it uniform} bounds on the kernel density of the semigroup $P_t$ with respect to the reference measure, see \cite{bcs,cks,coulhon,davies,va-sc-co} among many works on this topic. 

The Nash inequalities \eqref{Nash2} do not give the optimal constant $C$ in~\eqref{ineq:ultrac}. The optimal contractive bounds $\|P_tf\|_q\leq C_{p,q,n} (t) \|f\|_p$ for the classical heat equation in $\dR^n$  can be obtained by the Euclidean logarithmic Sobolev inequality (see~\cite{bakrystflour,ledoux2000}), but the Nash inequality is the easiest and the most intuitive way to get ultracontractive bounds such as~\eqref{ineq:ultrac}. 

\medskip

Inequalities~\eqref{Nash2} have been studied by F.-Y. Wang in~\cite{wang2000} as part of a more general family of inequalities, called  Super-Poincar\'e inequalities, of the form
\begin{equation}
\label{eq-nana}
\| f\|_2 ^2\leq a \, \cE(f,f) + b(a) \|f\|_1^2
\end{equation}
for $a > a_0$, where $b$ is a nonnegative function. Optimising in~\eqref{eq-nana} over the parameter $a$ leads to 
$$
\frac{\| f\|_2^2}{\|f\|_1^2 } \leq \psi \PAR{\frac{\cE(f,f) }{\|f\|_1^2}}
$$
where $\psi (x) = \displaystyle \inf_a \BRA{ a x + b(a) }$ is an increasing  concave function, or equivalently
\beq \label{Nash5}
\phi\PAR{\frac{\| f\|_2^2}{\|f\|_1^2 } }\leq \frac{\cE(f,f) }{\|f\|_1^2}
\eeq
  for an increasing convex function $\phi$. Then, following the argument leading to~\eqref{eq-ultra0}, it implies the ultracontractive bound
 \beq\label{ultrawang}
  \|P_t f\|_\infty\leq {U^{-1}(t)}\|f\|_1,
 \eeq
for all $t>0$,  where $U(t)=\int_t^\infty1/\phi(x)dx$ is well defined under adequate assumptions on $\phi$ (see \cite{coulhon}).

 The generalized Nash inequalities~\eqref{eq-nana} are also a powerful tool to obtain spectral properties of the generator defining the Dirichlet form (see~\cite{wang2000}); 
 in particular they imply that its essential spectrum is empty. When the reference measure has finite mass, they also provide additional properties of the measure in the fields of concentration, asymptotic behavior and isoperimetry, as in \cite{bcr07}. They belong to the large family of functional inequalities such as  the Logarithmic Sobolev and the Poincar\'e inequalities, and have been studied in many recent works such as \cite{maheux,wang-05}.
  
\bigskip

This work is devoted to a more general situation in which the semigroup is not ultracontractive, so that one cannot expect uniform bounds on its kernel density, as in~\eqref{ptunif}. For instance the Ornstein-Uhlenbeck on $\R^n$, which is probably the most studied semigroup on $\R^n$, beyond the classical heat semigroup, is not ultracontractive; in fact, according to a famous result by E.~Nelson, it  is only hypercontractive (see~\cite{logsob} for example). Observe, according to  the celebrated theorem of L.~Gross~\cite{gross75}, that the corresponding hypercontractive bounds are equivalent to a logarithmic Sobolev inequality for the Gaussian measure (which is weaker than the Sobolev inequality~\eqref{Nash3}). Of course the Ornstein-Uhlenbeck kernel is explicit, so it is useless to get any estimate on it, but, for many other examples, pointwise  estimates on the kernels are an interesting and not so easy issue. There is a very large literature on this problem, see \cite{davies} and the references therein.

Non-uniform estimates on the density of the heat kernel may provide useful information on the semigroup. For example, let us consider a symmetric semigroup $(P_t)_{t\geq0}$ which may be represented by a density $p_t(x,y)$ with respect to an invariant measure $\mu$, that is, such that
$$
P_tf(x)= \int_E f(y) p_t(x,y) d\mu(y)
$$
for all $x$ and $t>0$. Then the operator $P_t$ is in the trace class and therefore has a discrete spectrum as soon as $p_t(x,x) \in \cL^1(\mu)$~;  moreover estimates on the spectrum can be obtained as detailed below.

In the general situation when the kernel density should not be uniformly bounded, 
the classical Nash inequality~\eqref{Nash} is not adapted, and the main idea of this work is to use  the generalized  Nash inequality~\eqref{Nash5},   modified with a weight depending on the expected estimate. Depending on the generator of the heat kernel  and the reference measure considered in the $\mathcal L^p$ norms, we shall look for a positive function $V$ and an increasing and convex function $\phi$ such that
\begin{equation}
\label{eq-wnana}
\phi\PAR{\frac{\| f\|_2^2}{\|fV\|_1^2 } }\leq \frac{\cE(f,f) }{\|fV\|_1^2}
\end{equation}
for all $f$. Such an inequality will be called a {\it weighted Nash inequality}. We shall look for weight functions  $V$ satisfying the subharmonic condition $L V \leq c \, V$ where $L$ is the infinitesimal generator of the semigroup~; this assumption is very close (but easier to satisfy) to the condition on Lyapunov functions recently used by the first author, F. Barthe, P. Cattiaux and A. Guillin in \cite{guillin-al-08-1,guillin-al-08-2} to prove functional inequalities such as the Poincar\'e and super-Poincar\'e inequalities.  Here is a key difference between our approach and theirs~:  the Lyapunov functions used in the present work explicitly appear in the functional inequalities themselves, whereas in the works mentioned above they are only a tool
to get the sought functional inequalities but they do not explicitly appear in the final estimates~: they are used like a catalyst to derive them. We will prove that the weighted Nash inequality~\eqref{eq-wnana} and the subharmonic condition on the weight function $V$ imply the {\it non-uniform} estimate 
$$
p_t(x,y)\leq K(t, \phi, c)V(x)V(y)
$$
of the heat kernel, for a positive function $K$.

\section{Framework and outline of the work}

This work is devoted to properties of symmetric Markov semigroups $(P_t)_{t\geq0}$. On a given measure space $(E, \cE, \mu)$, a symmetric Markov semigroup is a family of positivity preserving operators acting on bounded measurable functions, which preserve  constant functions, and are moreover symmetric in $\cL^2(\mu)$.  In the main application of section~\ref{sec-app}, the measure $\mu$ will be a probability measure, but it could also be a measure with infinite mass.  The operators $P_t$ are contractions  in $\cL^1(\mu)$ and $\cL^\infty(\mu)$, so are contractions in any $\cL^p(\mu)$ with $1 \leq p \leq \infty$. The semi-group property consists in the identity $P_t\circ P_s = P_{t+s}$ for any $s$ and $t$ in $\R_+$, together with a continuity assumption at $t=0$, for example here that  for any $f\in \cL^2(\mu)$, $P_t f$ converges to $f$ in $\cL^2(\mu)$ when $t$ converges to $0$. We shall assume that, for all $t$, $P_t$ has a kernel, which is the case when $E$ is a Polish space. 

\medskip

Symmetric Markov semigroups naturally appear as the laws of Markov processes $(X_t)_{t\geq0}$ on $E$ which are reversible in time:  for example in the case  when $\mu$ is a probability measure, this means that for any $T>0$, the law of the process $(X_t, 0\leq t \leq T)$ when the law of $X_0$ is $\mu$  is the same as the law of the process $(X_{T-t}, 0\leq t\leq T)$. 

They also naturally  appear when solving a heat equation 
$$
\partial_t u= Lu, ~ u(x,0)= f(x)
$$
on $E\times [0, \infty)$; here $L$ is a  (unbounded) self-adjoint operator satisfying the maximum principle and $L1=0$, for example a second order differential sub-elliptic operator with no $0$-order term on an open set on $\R^n$ or a manifold; in this case, and under mild hypotheses, the solution may be represented as 
$$
  u(x,t) = P_t f (x).
  $$
  
  \medskip
  
 By the Hille-Yosida theory,   the operator $P_t$  has a derivative $L$ at $t=0$ which is defined in a domain dense in $\mathcal L^2(\mu)$. Moreover $P_t= \exp(tL)$ and $L$ is self-adjoint since $P_t$ is symmetric, see~\cite{yosida} for instance. Also $P_t$ is a contraction in $\cL^2(\mu)$, so that the spectrum of $L$ lies in $(-\infty, 0]$. 

\medskip

Under our assumptions, for all $t>0$ the operator $P_t$ will be represented by a kernel density $p_t(x,y)$ with respect to the reference measure $\mu$, in the sense that there exists a nonnegative symmetric function $p_t$ on $E \times E$ such that
 $$
 P_t f (x)= \int_E f(y) \, p_t(x,y) \, d\mu(y)
 $$
 for $\mu$ almost every $x$ in $E$.
Then the semigroup property $P_t\circ P_s = P_{t+s}$ may be translated into the celebrated Chapman-Kolmogorov equation
 $$
 \int_E p_t(x,y) p_s(y, z) d\mu(y)= p_{t+s}(x,z)
 $$
for $\mu \otimes \mu$ almost every $(x,z)$ in $E \times E$.

 Moreover, as soon as the  kernel density $p_t(x,y)$ is in $\cL^2(\mu\otimes \mu)$, the operator $P_t$ is Hilbert-Schmidt on $\mathcal L^2(\mu)$ (see~\cite{kolmogorov} for instance). In particular $P_t$ has a discrete spectrum  $(\mu_n(t))_{n\in\N}$, associated to a sequence of orthonormal eigenfunctions $(e_n)_{n\in\N}$ in $\cL^2(\mu)$. In this case 
 $$
 p_t(x,y) = \sum_n \mu_n(t) e_n(x)e_n(y)
 $$
  and the series converges since 
 \begin{equation}
 \label{eq-trace1}
 \sum_n \mu_n(t)^2 = \int_{E\times E} p_t(x,y)^2 d\mu(x)d\mu(y) < + \infty.
 \end{equation}
Moreover
$$
\int_{E\times E} p_t(x,y)^2 d\mu(x)d\mu(y) = \int_E p_{2t}(x,x) d\mu(x)
$$
so that $P_{2t}$ is in the trace class.  Of course such estimates can be established only for $t>0$.

  Since $P_t= \exp(tL)$ this just shows that $L$ itself has a discrete spectrum $(-\lambda_n)_{n\in \N}$ with $\lambda_n \geq 0$ and $\lambda_0=0$, such that $\mu_n(t) = e^{- \lambda_n t}$.
We see in the estimate~\eqref{eq-trace1} how a control on $p_t(x,x)$ or $p_t(x,y)$ may lead to a control on the spectrum $(\mu_n(t))_{n \in \N}$ of $P_t$, hence on the spectrum $(\lambda_n)_{n \in \N}$ of $L$.

 In general, as explained above, it is not easy to get the existence of the density $p_t(x,y)$ and such a control on it. The classical situation in which $P_t$ is Hilbert-Schmidt is when $\mu$ has finite mass and $p_t$ is bounded. For example, under the Nash inequality~\eqref{Nash} or~\eqref{Nash5}, then according to the ultracontractive bound~\eqref{ultrawang} the operator   $P_t$ is bounded from $\cL^1(\mu)$ into $\cL^{\infty}(\mu)$ with norm $C_t$. In this case $P_t$ may be represented by a kernel density $p_t$ which is $\mu \otimes \mu$ almost surely bounded by the same constant $C_t$ under a mild assumption on $(E, \cE, \mu)$ (for instance if $\cE$ is  generated by a countable family, up to zero measure sets, see~\cite[Lemma 4.3]{bakrystflour}):  spaces $(E, \cE, \mu)$ for which this holds will be called {\it nice measure spaces}.  They include Polish spaces  on which Markov semigroups  can be represented by a kernel. 
 
 \bigskip
 
This work is devoted to the case of non ultracontractive semigroups, that is, of non bounded kernel densities. 
We shall replace the Nash inequality by the weighted Nash inequality~\eqref{eq-wnana} with a weight $V$ such that $L V \leq c V$ to obtain the existence of a density $p_t$ which satisfies 
\begin{equation}
\label{eq-derder}
p_t(x,y) \leq K(t,\phi, c) \, V(x) \, V(y),
\end{equation}
see Proposition~\ref{prop:WUltraContr}, Theorem~\ref{thm:BornesPt} and Corollary~\ref{coro-sect2}. 

\medskip

In section~\ref{sec-univ} we give a simple illustration of this method, see Theorem~\ref{thm:bornesUniv}. There we deduce the following universal bound on $\dR^n$ from the classical Nash inequality~\eqref{Nash}~: if the invariant measure $\mu$, not necessarily  finite, has a positive density $\rho$, then    
$$
 \vert\vert f\vert\vert_{2}^{2+\frac{4}{n}}
\leq C_n\,\vert\vert f V\vert\vert_{1}^{\frac{4}{n}}\, \PAR{
\cE(f,f)+
\int_{\R^n} \frac{LV}{V} f^2\, d\mu},
$$
where $V=\rho^{-1/2}$. This leads to a weighted  Nash inequality if moreover  $L V \leq c V$, whence to bounds such as~\eqref{eq-derder}.

\medskip

  A case study  of symmetric semigroups on $\dR$ consists in the Sturm-Liouville operators : given a  probability measure $\mu$ with  smooth and positive density $\rho$ with respect to the Lebesgue measure, the Sturm-Liouville operator  
 $$
 L f =  f'' + \log (\rho)'f'
 $$  
 defined on smooth functions leads to a symmetric Markov semigroup  in $\mathcal L^2(\mu)$. Depending on $\rho$, this family shows all possible behaviours. 
 The main example studied in this article concerns the probability measures 
  $$
d\mu_a(x)= \rho_a(x)dx=C_a e^{-|x|^a}dx
 $$
on $\dR$ and their associated Markov semigroup $(P_t)_{t\geq0}$; here $a>0$ and $C_a$ is a normalization constant. 

 If $a>2$ then the  semigroup is ultracontractive and the density with respect to the measure  $\mu_a$ is uniformly bounded (see~\cite{kkr} for the proof, among more general examples). In the limit Gaussian case when $a=2$ then the semigroup is the well known Ornstein-Uhlenbeck semigroup (up to normalization),  
which is not ultracontractive any more  but only hypercontractive.  It means that for $t>0$, $P_t$ maps $\mathcal L^2(\mu_a)$ into some $\mathcal L^{q(t)}(\mu_a)$, where $2<q(t)<\infty$ :  this is Nelson's Theorem. Observe that in this case one explicitly  knows the density $p_t(x,y)$ and the spectrum $\lambda_n= n$, and that $P_t$ is Hilbert-Schmidt. 
  
Now, if $1<a<2$ the semigroup $P_t$ is not hypercontractive anymore since the measure $\mu_a$ does not satisfy a logarithmic Sobolev inequality anymore. In fact, as shown in~\cite{bcr1}, $P_t$ with $t>0$ satisfies Orlicz hypercontractivity : it maps $\mathcal L^2(\mu_a)$ into a Orlicz space  slightly smaller than    $\mathcal L^2(\mu_a)$. This functional regularity does not bring any explicit upper bound on the kernel density $p_t$.   
  
As a simple illustration of our general method,  we shall prove  that for all real $\beta$ there exists $\theta >0$ such that the density $p_t(x,y)$ satisfies the explicit upper bound 
  $$
p_t(x,y)\leq  C(a,\beta)\frac{e^{ct}}{t^{\theta }}\frac{\rho_a^{-1/2}(x)\rho_a^{-1/2}(y)}{(1+|x|^2)^{\beta}(1+|y|^2)^{\beta}}.
  $$
 For $\beta>1/2$, this estimate is in $\mathcal L^2(\mu_a)$, so that the operator $P_t$ is Hilbert-Schmidt : to our knowledge this is a new result.  
 In the other limit case, when  $a=1$, such estimate  can not hold anymore : indeed the spectrum of $-L$ does not only have a discrete part  but lies in $\{0\}\cup (\lambda_0, \infty)$, with $\lambda_0>0$ (see \cite{wangbook}).  Let us note that studying the measures $\mu_a$ for $a\in(1,2)$ is a current active domain in functional analysis. These measures represent a large class of log-concave measures:  they are not log-concave enough to satisfy a logarithmic Sobolev inequality, but some of their properties, as the concentration for instance, are similar of the standard Gaussian measure, one can see \cite{bcr1,bcr07,ggm06,latala} for example.

  \medskip
  
  The method used here to get the weighted Nash inequalities on the real line will be quite close to the method introduced by B. Muckenhoupt in~\cite{muckenhoupt} and generalized later by S. Bobkov and F. G\"otze in~\cite{bobkov-gotze} to characterize measures which satisfy Poincar\'e or logarithmic Sobolev inequalities in the real line. We shall not try here to get the same kind of if and only if results, since there are too many parameters to control (the weight function $V$, the rate function $\Phi$ and so on). 
  
 \medskip
 
 We shall not either try  to extend our results to the most general setting, for example Riemannian manifolds, which would require a more precise analysis of the Laplacian of the distance function, and therefore lower bounds on the Ricci curvature.  Instead we prefer to concentrate on some key one-dimensional models  to show the easiness and the efficiency of the methods presented here. Moreover, as  usual when using Lyapunov functions, constants obtained in these estimates are far from optimal and that is why we  only  focus on the overall  behavior of the estimates but not try to make the constants finer.

 \medskip

The plan of the article is the following. In the next section we explain the abstract result :  how a weighted Nash inequality coupled to a Lyapunov function implies a non-uniform estimate of the kernel density. In section~\ref{sec-univ} we prove a universal weighted Nash inequality. 
In section~\ref{sec-app} we finally apply the method of section~\ref{sec-prems} to the measures  $\mu_a$ defined above for $a\in(1,2)$.

\medskip

\noindent {\bf Notation :} In the whole  article, $\|\cdot\|_p$ stands for  the $\cL^p$ norm  with respect to the  measure $\mu$. The measure $\mu$ could change, depending on the context,  but it should be always clear.
 \section{The abstract  result}
 \label{sec-prems}
 
 In this section we present a simple method to obtain the existence and explicit and non-uniform bounds on  Markov semigroup kernel densities. 
 
In the classical ultracontractive  case the upper bound on the kernel density $q^2$ of $Q\circ Q$ follows from
$$
 \|Qf\|_2 \leq \| f\|_1\Leftrightarrow \|Q\!\circ\! Q\,f\|_{\infty} \leq \| f\|_1\Leftrightarrow |q^2(x,y)| \leq 1.
$$
We extend this property to non-uniform estimates.

 \bprop\label{prop:WUltraContr} Let $(E, \cE, \mu)$ be a nice measure space, $Q$ a symmetric bounded operator on $\cL^2(\mu)$ and  $V$ a positive measurable function on $E$. Then the two assertions are equivalent :
 \begin{enumerate}[(i)]
 \item The operator $Q$ satisfies 
$$
 \|Qf\|_2 \leq \| fV\|_1
$$
for all $f\in \cL^2(\mu)$ ;
\item  The operator $Q^2= Q\circ Q$ may be represented by a kernel density $q^2(x,y)$ with respect to $\mu$ which satisfies
$$
 |q^2(x,y)| \leq V(x)V(y)
$$
 for $\mu \otimes \mu$ almost every $(x,y)$ in $E \times E$. 
 \end{enumerate}
If moreover the function $V$ is in $\cL^2(\mu)$, then $Q$ is Hilbert-Schmidt, and therefore has a discrete spectrum $(\mu_n)_{n\in\N}$ such that,
 $$
 \sum_n \mu_n^2 \leq \int V^2 d\mu.
 $$
 \eprop
  \bpf
 Let us assume $(i)$ and let us consider the operator $Q_1= \frac{1}{V}Q V$, that is, defined by 
 $$
 Q_1 f= \frac{1}{V} Q(fV).
 $$ 
 By hypothesis, $Q_1$  is a contraction from $\cL^1(\nu)$ into $\cL^2(\nu)$ where $d\nu = V^2d\mu$. Moreover it is symmetric with respect to the measure $\nu$ since so is $Q$ with respect to $\mu$, so by duality it is also a contraction from
  $\cL^2(\nu)$ into $\cL^\infty(\nu)$, and by composition the operator $Q^2_1=Q_1\!\circ\! Q_1$ is a contraction from $\cL^1(\nu)$ into $\cL^\infty(\nu)$.  

 This implies that $Q_1^2$ may be represented by a kernel density $q_1^2(x,y)$   in the space $\cL^2(\nu)$   which satisfies  $|q_1^2(x,y)|\leq 1$ for $\nu \otimes \nu$ almost every $(x,y)$ in $E \times E$ (see \cite[Lemme 4.3]{bakrystflour} for instance). On the other hand, 
 $$
 q_1^2(x,y) \, V(x) \, V(y) = q^2(x,y)
 $$
for $\mu \otimes \mu$ every $(x,y),$ noting that $V$ is positive. This implies $(ii)$. 
 
 \medskip
 
 Conversely, if $f \in \cL^2(\mu)$, then, by symmetry of $Q$,
 $$
 \|Qf\|^2_2= \int fQ^2 f  \, d\mu = \int q^2(x,y) \, f(x) \, f(y) \, d(\mu \otimes \mu) (x,y) \leq \PAR{\int \vert f \vert V \, d\mu }^2,
 $$
 which proves $(i)$.
 
\medskip 
  
  If now $V\in \cL^2(\mu)$, then the kernel  $q^2(x,x)$ is integrable on $E $ with respect to $\mu$, which just means that $Q$ is Hilbert-Schmidt.
 \epf

\begin{ex}
The first and explicit example is the classical Ornstein-Uhlenbeck semigroup in $\dR^n$, with generator $L = \Delta - x \cdot \nabla$ : in a probabilistic form it is given by the Mehler formula
 $$
 P_t f (x) = \E\PAR{f(e^{-t}x+ \sqrt{1-e^{-2t}}Y)},
 $$
  where $Y$ is a standard Gaussian variable with law $\gamma$. It admits a kernel density with respect to the Gaussian measure, given by
   $$
 p_t(x,y)= (1-e^{-2t})^{-n/2}\exp\SBRA{-\frac{1}{2(1-e^{-2t})}(|y|^2e^{-2t} -2 \, x \cdot ye^{-t}+ |x|^2e^{-2t})}
 $$ 
 for all $x,y\in\dR^n$ and $t>0$. In particular
\begin{equation}\label{boundNt}
p_{2t}(x,y)\leq p_{2t} (x,x)^{1/2} p_{2t} (y,y)^{1/2} = (1-e^{-4t})^{-n/2} \exp\PAR{\frac{|x|^2}{1+e^{2t}}} \exp\PAR{\frac{|y|^2}{1+e^{2t}}}
\end{equation}
 by the Cauchy-Schwarz inequality, with equality if $x=y$. Hence, by Proposition~\ref{prop:WUltraContr},
 $$
\NRM{P_t{f}}_{L^{2}(d\gamma)}\leq \|fV_t\|_{L^1(d\gamma)}
$$
where 
$$
V_t(y)=(1-e^{-4t})^{-n/4}\exp \PAR{\frac{|y|^2}{2(1+ e^{2t})}}.
$$
This bound has been obtained in a more general context in~\cite{bbg1}, where it is shown to be optimal, being an equality for square-exponential functions $f$.
\end{ex}

By Proposition~\ref{prop:WUltraContr} we are now brought to prove bounds such as $(i).$

When the operator $Q$ is a Markov semigroup $P_t$ with a kernel $p_t$, evaluated at time $t$, then one may obtain such bounds through functional inequalities that we describe here. We shall mainly be concerned with the case when $\mu$ is a probability measure, although much of what follows could be extended to the case when $\mu$ has infinite mass.

Let  $(P_t)_{t\geq0}$ be  a symmetric Markov semigroup on $E$ with generator $L$ and associated Dirichlet form
$$
\cE_{\mu}(f,f) = -\int f Lf d\mu.
$$
This quadratic form can be defined  on a larger subspace than the domain of $L$, which is called the domain of the Dirichlet form.

\medskip

Bounds such as $\Vert P_t f \Vert_2 \leq K(t) \Vert f V \Vert_1$ will be obtained by means of weighted Nash inequalities and Lyapunov functions, that we now define. 

\bdefi 
Let $V$ be a positive function on $E$, $M$ be a nonnegative real number and  $\phi$ be a positive function defined on $(M, \infty)$ with $\phi(x)/x$ non decreasing.

The Dirichlet form $\cE_{\mu}$ satisfies a weighted Nash inequality  with weight $V$ and rate function~$\phi$ if
\beq 
\label{ineg:NashPhi} 
\phi \PAR{ \frac{\|f\|_2^2}{\|f V\|_1^2} } \leq \frac{\cE_\mu(f,f)}{\|f V\|_1^2}
\eeq
for all functions $f$ in the domain of the Dirichlet form such that $\displaystyle \|f\|_2^2 > M \, \|f V\|_1^2.$ 
\edefi

As recalled in the introduction, the fundamental two examples are the classical Nash inequality~\eqref{Nash} for the Lebesgue measure, with $\phi(x) = C x^{1+2/n}$ and $M=0$, ($M>b^{n/2}$ for the generalized inequality~\eqref{Nash2}) and those~\eqref{Nash5} given by Super-Poincar\'e inequalities, with $\phi$ the inverse of $\displaystyle \inf_a\BRA{ a x + b(a)}$ and $M=0$. They all have weights $V=1$, and in the following we shall be concerned with Nash inequalities with a general positive weight $V.$

\bdefi 
\label{def-2}
A Lyapunov function is a positive function $V$ on $E$ in the domain of the generator $L$ such that
\beq \label{ineq:subharm} 
LV \leq cV
\eeq 
for a real constant $c$, called the Lyapunov constant.
\edefi

It is not really necessary for $V$ to be in the $\cL^2$-domain of $L$, but for simplicity we restrict to this situation, which will be the situation in our examples below.

\brmq
\label{rem-33}
In our context the Lyapunov constant $c$ will be nonnegative. Negative Lyapunov constants can also be considered, but by adding an extra term :  for instance the authors in~\cite{guillin-al-08-1,guillin-al-08-2} consider Lyapunov functions $V$ such that $LV\leq -\gamma V+\I_K$ where $\gamma>0$, $V\geq 1$ and $K$ is a compact set.  These Lyapunov functions are a powerful tool to obtain  rates of the long time behavior of the Markov semigroup, for example, through the obtention of Poincar\'e or more generally weak Poincar\'e inequalities.

As mentioned in the introduction, Lyapunov functions defined as in our definition~\ref{def-2} with $c \geq 0$ are introduced to obtain smoothing properties of the Markov semigroup for a fixed time $t>0$. 
 \ermq

 When $\mu$ has finite mass,  one can also observe that the restriction $V\geq0$ in~\eqref{ineq:subharm}  could be replaced by  $V\geq1$ when  $c\geq 0$, since one may always change $V$ into $V+1$. This will be the case in the main application given in section~\ref{sec-app}.

Then, one has the following.
\bthm[Wang]
\label{thm:BornesPt}
Let $(P_t)_{t\geq 0}$ be a Markov semigroup on $E$ with generator $L$ symmetric in $\cL^2(\mu)$. 

Assume that there exists a Lyapunov function $V$ in $\mathcal L^2(\mu)$ with Lyapunov constant $c\geq 0$, and that the Dirichlet form associated to $L$ satisfies a weighted Nash inequality with weight $V$ and rate function $\phi$ on $(M, +\infty)$ such that 
\begin{equation}
\label{eq-hypom}
\int^\infty \frac{1}{\phi(x)} dx < \infty.
\end{equation}
Then 
$$
\|P_tf\|_2 \leq K(2t) \,  e^{ct} \Vert f V \Vert_1
$$
 for all $t >0$ and all functions $f\in \cL^2(\mu)$; here  the function $K$  is defined by
$$
K(x)=
\left\{
\begin{array}{ll}
\sqrt{U^{-1}(x)}\,\,&{\rm if} \,\, 0 < x < U(M),\\
\sqrt{M}\,\,&{\rm if}\,\,x \geq U(M)
\end{array}
\right. 
$$
where $U$ denotes the (decreasing) function defined on $(M, +\infty)$ by
$$
U(x)= \int_x^\infty \frac{1}{\phi (u)} du.
$$  
\ethm

\brmq
After completing this work, we learnt from F.-Y. Wang that he obtained this  result  under weighted Super-Poincar\'e inequalities in~\cite[Theorem~3.3]{wang02}. We state  and prove it in our context to show that our method is simple and self contained.  
\ermq

 Here the measure $\mu$ need not be a probability measure and  may have infinite mass and, in the case when   $U(M) = +\infty$, then $K$ is just defined by the first line.

Observe also that if $M=0$ then we can take any real parameter~$c$, as one can see from the proof. 

\brmq
As mentionned in Remark~\ref{rem-33}, we are not mainly concerned with the long time behaviour of the Markov semigroup, though in some cases a weighted Nash inequality may reveal adapted: for instance, in the case when $c=0, M=0$ and $U(M)=0$, then Theorem~\ref{thm:BornesPt} ensures that $P_t f$ converges to $0$ in $\cL^2(\mu)$ for all $f \in \cL^2(\mu)$ with finite $\Vert f V \Vert_1$; observe that in this case $\mu$ has necessarily infinite mass. If $\mu$ is a probability measure, then we expect $P_t f$ to converge to $\int fd\mu$, which is a priori nonzero, so the rate $K(2t)e^{ct}$ can not  converge to 0.

On the contrary weighted Nash inequalities are adapted to get estimates on the small time behavior : Theorem~\ref{thm:BornesPt} gives a bound on $\Vert P_t f \Vert_2$ for $t>0$ which depends on $f$ only in terms of a weighted $\cL^1$ norm, which is an illustration of the gain of integrability induced by the semigroup. Observe that the coefficient $K(2t)$ tends to $+ \infty$ as $t$ goes to $0$.
\ermq

By Proposition~\ref{prop:WUltraContr} this leads to the following bounds on the kernels:

\bcor\label{coro-sect2}
If the Markov semigroup $(P_t)_{t\geq0}$  satisfies the assumptions of Theorem~\ref{thm:BornesPt} above, then $P_t$ has a density $p_t$  with respect to $\mu$ which satisfies 
$$
 p_{2t}(x,y)\leq K(2t)^2e^{2ct}V(x)V(y),
$$
for all $t>0$ and $\mu\otimes\mu$ almost every $(x,y)\in E \times E$. 

Moreover $P_t$  is Hilbert-Schmidt for all $t>0$, and therefore has a discrete spectrum $(\mu_n(t))_{n\in\N}$ such that
 $$
 \sum_n \mu_n(t)^2 \leq K(2t)^2 e^{2ct} \int V^2 d\mu.
 $$
\ecor

\noindent{\it Proof of Theorem~\ref{thm:BornesPt}. ---}
Let $f$ be given in $\cL^2(\mu)$.
With no loss of generality we can assume that $f >0$ by writing the argument for $\vert f \vert + \varepsilon$, and letting $\varepsilon$ go to $0$ and using the bound $\vert P_t f \vert \leq P_t \vert f \vert.$

First notice that the map $G(t)= \int V P_t f d\mu$ has derivative
$$
G'(t)= \int VLP_t f d\mu = \int LV P_t f d\mu\leq c \, G(t),
$$
so that 
\beq 
\label{ineg:G}
\int V P_t f d\mu \leq e^{ct} \int V f d\mu.
\eeq

Then, given $0 \leq t \leq T$ fixed, consider the function 
$$
R(s)= \frac{\|P_sf\|_2^2}{\PAR{e^{ct}\int fVd\mu}^2}
$$
on $[0,t].$ Then
\begin{equation}
\label{eq-derive}
\frac{-R'(s)}{2}= \frac{\cE_{\mu}(P_sf, P_s f)}{\PAR{e^{ct}\int fVd\mu}^2}=\frac{\cE_{\mu}(P_sf, P_s f)}{\PAR{\int P_s fVd\mu}^2}\PAR{\frac{\int P_sfVd\mu}{e^{ct}\int fVd\mu}}^2.
\end{equation}
In particular $R$ is decreasing. Moreover, if there exists $s\in [0,t]$ such that $R(s) \leq M$, then $R(t) \leq R(s) \leq M$, which yields the result. Hence we now assume that  $R (s) \geq M$ on $[0,t]$. Then, by ~\eqref{ineg:G},
$$
\frac{\|P_sf\|_2^2}{\PAR{\int P_s fVd\mu}^2} 
=
\frac{\|P_sf\|_2^2}{\PAR{e^{ct} \int f Vd\mu}^2} \frac{\PAR{e^{ct} \int f Vd\mu}^2}{\PAR{\int P_s fVd\mu}^2}
=
R(s) e^{2c(t-s)} \frac{\PAR{e^{cs} \int f Vd\mu}^2}{\PAR{\int P_s fVd\mu}^2}
\geq M
$$
for $c \geq 0.$ 

Hence, by applying the weighted Nash inequality to $P_s f$,~\eqref{eq-derive} gives 
$$
\frac{-R'(s)}{2}\geq \phi\PAR{\frac{\|P_sf\|_2^2}{\PAR{\int P_s fVd\mu}^2}} \PAR{\frac{\int P_sfVd\mu}{e^{ct}\int fVd\mu}}^2.
$$
Moreover
$$
 \phi\PAR{\frac{\|P_sf\|_2^2}{\PAR{\int P_s fVd\mu}^2}}\geq \phi\PAR{\frac{\|P_sf\|_2^2}{\PAR{e^{ct}\int fVd\mu}^2}}\PAR{\frac{\int P_sfVd\mu}{e^{ct}\int fVd\mu}}^2
$$
from the inequality \eqref{ineg:G} and the fact that $\phi(x)/x$ is non decreasing, so that
$$
\frac{-R'(s)}{2}\geq \phi(R(s)).
$$
In turn this may be seen as
$$
U(R(s))' \geq {2}
$$ 
which integrates into
$$
U(R(t))\geq U(R(0))+{2t}\geq {2t}.
$$
Since $U^{-1}$ is defined on $(0,U(M)]$ and is decreasing then we obtain the upper bound 
$$ 
R(t) \leq U^{-1}(2t)
$$
for all $t\leq U(M)/2$. For those $t\geq U(M)/2$ then we have $R(t) \leq M$. Combining all these estimates gives the result.  
{\hfill\rule{1.3mm}{2mm}\par} 

\brmq\label{L2weight}
In the main application of the weighted Nash inequality given in section~\ref{sec-app},  the weight function $V$ is in $\mathcal L^2(\mu)$. But formally, one does not need $V$ to be in $\cL^2(\mu)$ to get the result. This restriction is made here not only in view of Proposition \ref{prop:WUltraContr}. It is also made to ensure the integration by parts formula
$$
\int LP_s f V d\mu = \int P_s f LV d\mu
$$
which leads to~\eqref{ineg:G}, and automatically holds when $V$ is in $\cL^2(\mu)$ and in the domain of $L$. For those $V$ which increase too rapidly at infinity, then it may be false in general; it requires a more precise analysis of the semigroup $(P_t)_{t\geq0}$ and restricting to a large  subclass of functions in $\cL^2(\mu)$. 

Here are two fundamental examples in the two cases when $\mu$ has finite or infinite mass~:
\begin{itemize}
\item The Lebesgue measure on $\dR^n$  satisfies the classical Nash inequality~\eqref{Nash}, hence a weighted Nash inequality with weight $V=1$ and rate function $\phi(x)= C x^{1+2/n}$, for instance on the set $(0,+\infty)$. Then, by Theorem~\ref{thm:BornesPt} applied with $V=1$ and $c=0$, one recovers the well known contraction property of the classical heat kernel on $\dR^n$, 
$$
\|P_tf\|_2\leq \Big(\frac{C}{t}\Big)^{n/4}\|f\|_1,
$$
for all $t>0$ and for the non optimal constant $C = n/4$ instead of $1/(8 \pi)$ (see~\cite{ledoux2000} for instance). In this case, we only have to consider functions $f\in\mathcal L^1(\mu)$ and in the domain of the  Dirichlet form $\int \vert \nabla f \vert^2 d\mu.$   The main tool to get optimal bounds in any $\mathcal L^p(\mu)$ for $p\geq 1$ for the classical heat kernel on $\dR^n$ is the Euclidean logarithmic Sobolev inequality as explained for instance in~\cite{bakrystflour} or~\cite{ledoux2000}.   

\item The second example concerns the Sturm-Liouville operator $Lf = f'' + (\log \rho)' f'$ on $\R$, associated with the measure $d\mu= \rho(x)dx$. Here it would be enough to know that $(\log \rho)''$ is bounded from above and that $V\rho'$ and $V'\rho$  go to $0$ at infinity. Indeed, in this situation, it is enough for smooth functions $f$ and $g$  that $f'g\rho$ and $f g' \rho$ go to $0$ at infinity to ensure, through integration by parts, that
$$\int Lf g d\mu= - \int f' g' d\mu = \int fLg d\mu.$$
When $(\log \rho)''$ is bounded from above, the semi-group satisfies a $CD(a, \infty)$ inequality; hence, as soon as $f$ is bounded, then so is $(P_t f)'$ when $t>0$ (see~\cite[Remark 5.4.2]{logsob}). Hence in this case we may work with the space of bounded functions to get the result. 

Examples will be studied in sections~\ref{sec-univ} and~\ref{sec-app}.  
\end{itemize}
\ermq

Theorem \ref{thm:BornesPt} has the following converse:

\bthm\label{thm:RcqueNash} 
Let $\mu$ be a measure on $E$ and let $(P_t)_{t\geq 0}$ be a Markov semigroup on $E$ with generator $L$ symmetric in $\cL^2(\mu)$.  

If there exists a positive  function $V$ and a positive function $K$ defined on $(0,\infty)$ such that
$$
\|P_tf\|_2 \leq K(t)\| fV\|_1
$$ 
for all $t > 0$, then the weighted Nash inequality \eqref{ineg:NashPhi}  holds with the same function $V$, $M=0$ and function 
$$
\phi(x)= \sup_{t>0} \frac{x}{2t}\log\frac{x}{K(t)^2}, \quad x \geq 0.
$$
\ethm

Here again $\mu$ need not be a probability measure.

\brmq For instance, by Theorem~\ref{thm:BornesPt}, if we assume a Nash inequality with $\phi(x)= Cx^r$ for large $x$,  with $r>1$, then we obtain a bound such as $\|P_tf\|_2 \leq K(t)\| fV\|_1$ with $K(t)= C't^{1/2(1-r)}$ for small $t$.
 
 Conversely,  if we assume such a bound with such a $K$, then, by the converse Theorem \ref{thm:RcqueNash}, we obtain a Nash inequality with function $\phi (x) = C'' x^{r}$ for large $x$. Therefore, in this case and up to the values of the constants, we have a true quantitative equivalence between the Nash inequality and the bound on $\|P_t f\|_2$. 
\ermq

{\noindent{\it Proof of Theorem~\ref{thm:RcqueNash}. --- }}
It is based on the observation that the function 
$$
t\mapsto \log (\|P_tf\|_2^2)
$$
 is convex for any symmetric semigroup.  Indeed, if 
$h(t)= \| P_t f\|^2,$
then
 $h'(t)= 2 \int P_t f L(P_tf) d\mu$ and 
$h''(t)= 4 \int (LP_t f)^2 d\mu$; hence
$h'^2 \leq h h''$,  or equivalently $(\log h)'' \geq 0$. 

Therefore
$$
\log h(u) - \log h(0) \leq \frac{u}{t} \big[ \log h(t) - \log h(0) \big]
$$
for all $0 < u \leq t$, so that
 \begin{equation}
 \label{eq-conv10}
 h'(0)\leq \frac{h(0)}{t}\log \frac{h(t)}{h(0)}
 \end{equation}
by letting $u$ go to 0. 

Now, if moreover 
$$
h(t) \leq K(t)^2 \|fV\|_1^2,
$$
then~\eqref{eq-conv10} gives
$$
-2\frac{\cE(f,f)}{\|fV\|_1^2}\leq \frac{\|f\|_2^2}{\|fV\|_1^2}\frac{1}{t}  \log\PAR{\frac{K(t)^2\|fV\|_1^2}{\|f\|_2^2}}.
$$
This gives the claimed weighted Nash inequality.
{\hfill\rule{1.3mm}{2mm}\par}  
 
%%%%%%%%
\section{A universal weighted Nash inequality on $\dR^n$} 
\label{sec-univ}

Let $\rho$ be a positive smooth function on $\dR^n$.  We prove  a weighted Nash inequality   for   the operator  $Lf=\Delta f+\nabla \log \rho\cdot \nabla f$  with  the universal weight 
 $V={\rho}^{-1/2}$ and the measure $d\mu(x)=\rho(x)\,dx$. As usual, $\|\cdot\|_p$ stands for the $\mathcal L^p(\mu)$ norm and $(P_t)_{t \geq 0}$ is the semigroup with generator $L$.

\bthm\label{thm:bornesUniv} 
 In the above notation, the classical  Nash inequality~\eqref{Nash} is equivalent to
\begin{equation}
\label{nashmu}
 \vert\vert f\vert\vert_{2}^{2+\frac{4}{n}}
\leq C_n^{\frac{4}{n}}\,\vert\vert f V\vert\vert_{1}^{\frac{4}{n}}\, 
\left( 
\cE(f,f)+
\int_{\R^n} \frac{LV}{V} f^2\, d\mu
\right)
\end{equation}
for all smooth functions $f$ on $\dR^n$ with compact support. If moreover ${LV}\leq cV$ for $c\in\dR$ then 
$$
 \vert\vert f\vert\vert_{2}^{2+\frac{4}{n}}
\leq C_n^{\frac{4}{n}}\,\vert\vert f \, V \vert\vert_{1}^{\frac{4}{n}}\, 
\left( \cE(f,f) + c
\int_{\R^n}   f^2\, d\mu
\right)
$$
\ethm

\bpf
Let $g$ be a smooth function with compact support and let  $f=g\sqrt{\rho}$. Then 
$$
\int_{\R^n} \vert f\vert^2\, dx=  \vert\vert g\vert\vert_{2}^2,
$$
$$
\int_{\R^n} \vert f\vert\, dx= \int_{\R^n} \vert g\vert \sqrt{\rho}
\, dx= 
 \vert\vert g V\vert\vert_{1},
 $$
 and
 $$
 \int_{\R^n} \!\vert \nabla f\vert ^2\, dx
% =
%  \int_{\R^n} \!\ABS{\nabla \PAR{\frac{g}{V}}}^2\, dx
   = \int_{\R^n} \!\vert \nabla g\vert^2\, d\mu+
     \int_{\R^n} \! 2 \frac{g}{V} \nabla g.   \nabla  \frac{1}{V}\, dx
     +   \int_{\R^n}  
    g^2 
  \ABS{\nabla  \PAR{\frac{1}{V}}}^2\, dx.
 $$
By integration by part, the middle term is
$$
   \int_{\R^n}  \nabla (g^2). \PAR{  \frac{1}{V} \nabla  \frac{1}{V}} dx= -\int_{\R^n} g^2  \nabla\PAR{ \frac{1}{V} \nabla  \frac{1}{V}} dx  =\int_{\R^n} g^2 
 \left( \frac{\Delta V}{V} -3\frac{ \vert \nabla V\vert^2}{V^{2}}\right)d\mu,
$$
%The last term
%$$
%\int_{\R^n}  
%    g^2 
%  \vert\nabla  \frac{1}{V}\vert^2dx=
%\int_{\R^n} g^2 
%\frac{ \vert \nabla V\vert^2}{V^{2}}
% d\mu
%$$
so that
$$
 \int_{\R^n} \vert \nabla f\vert ^2\, dx
 = \int_{\R^n} \vert \nabla g\vert^2\, d\mu +
        \int_{\R^n} g^2 
 \left( 
  \frac{\Delta V}{V}
   -
   2\frac{ \vert \nabla V\vert^2}{V^{2}}\right)
  \, d\mu.
 $$
Moreover  
$$
\frac{LV}{V}= \frac{1}{V}\PAR{ \Delta V-2\,\nabla\log{{V}}\cdot\nabla V }=
  \frac{\Delta V}{V}
   -
   2\frac{ \vert \nabla V\vert^2}{V^{2}},
  $$
so
$$
 \int_{\R^n} \vert \nabla f\vert ^2\, dx=\cE(g,g)+\int_{\R^n} \frac{LV}{V} g^2\, d\mu.
$$
Hence the classical Nash inequality~\eqref{Nash} for $f$ is equivalent to~\eqref{nashmu} for $g$, which concludes the proof.
\epf

\bigskip

This type of transformation has been performed by F.-Y. Wang in~\cite{wang02} at the level of the Super-Poincar\'e inequality~\eqref{eq-nana}. From this the author estimates the kernel density of semigroups with infinite invariant measure. From Theorem~\ref{thm:bornesUniv} we now give estimates in the case of probability invariant measures. 

\begin{cor}\label{cor-VL1}
In the above notation, assume that $\mu$ is a probability measure and that $V\in\mathcal L^1(\mu)$ satisfies $LV\in\mathcal L^1(\mu)$ and $LV\leq cV $ with $c\geq0$. Assume moreover that the Hessian of $\log \rho$ is uniformly bounded from above on $\dR^n$ and that
$$
\sup_{\vert x \vert = r} \rho(x)^{1/2} \, r^{n-1} \to 0 \qquad \textrm{and} \qquad \sup_{\vert x \vert = r} \vert \nabla \rho(x) \vert \rho^{-1/2} \, r^{n-1} \to 0
$$
as $r$ tends to infinity. Then $P_t$ has a density $p_t$ which satisfies
 \begin{equation}\label{}
p_{2t}(x,y)\leq {\frac{d}{t^{n/2}}\, e^{2ct}}{}\,V(x)V(y)
\end{equation} 
for some $d>0$ and for all $x,y\in \R^n$, $t>0$.
\end{cor}

\bpf
We cannot directly apply Theorem~\ref{thm:BornesPt} since $V={\rho}^{-1/2}$ is never in ${\mathcal L}^2(\mu)$. The argument is exactly the same, but we have to justify the inequality $G'(t) \leq c \, G(t)$ where $G(t) = \int V \, P_t f \, d\mu$ for any smooth function $f$ with compact support. 

First of all $G'(t) = \int V L P_t f \, d\mu$ since $V \in \mathcal L^1(\mu)$ and $Lf$ is bounded.

Then we prove the integration by parts
$$
\int_{\dR^n} V L P_t f \, d\mu = \int_{\dR^n} L V P_t f \, d\mu.
$$
Let $r>0$, $B_r$ be the centered ball of $\dR^n$ with radius $r$ and $\vec{v}$ be its outward unit normal vector. Then, by two integrations by parts on $B_r$,
\begin{multline*}
\int_{B_r}  \! \! V L P_t f d\mu =  \int_{B_r}\! \! L V P_t f d\mu \\
- \int_{S^{n-\!1}} \! \! \! \!P_t f (r \omega) \, \nabla V (r \omega) \cdot \vec{v}\, \rho(r \omega) r^{n-1}\,d\omega + 
\int_{S^{n-\!1}}\! \!\! \! V (r \omega)\, \nabla P_t f (r \omega) \cdot \vec{v}\, \rho(r \omega) r^{n-1}\,d\omega.
\end{multline*}

But the Hessian of $\log\rho$ is uniformly bounded from above on $\dR^n$, say by the real number $\lambda$, so $L$ satisfies a $CD(- \lambda, \infty)$ curvature-dimension criterion. In particular (see~\cite{bakrystflour} for instance) it implies the uniform bound
$$
| \nabla P_t f | \leq e^{\lambda t} \, P_t | \nabla f | \leq e^{\lambda t} \, \Vert \nabla f \Vert_{\infty}.
$$
Then our assumptions on $\rho$ ensure that the last two terms tend to $0$ as $r$ tends to infinity, which justifies the integration by parts.
\epf

\brmq
The key point here is that $V={\rho}^{-1/2}$ is never in ${\mathcal L}^2(\mu)$, so this result does not ensure whether $P_t$ is  Hilbert-Schmidt or not. 

\ermq

We illustrate Corollary~\ref{cor-VL1} on the examples of Cauchy and exponential type measures. We have in mind the measure  $\exp(-|x|^a) dx$ in $\dR^n$ but for convenience we will study $\exp(-(1+|x|^2)^{a/2}) dx$ instead of $\exp(-|x|^a) dx$ which has the same behavior at infinity and has no singularity at $x=0$.

\bcor\label{appliCauchyGauss}
Let $\rho(x)=(1+\vert x\vert ^2) ^{-\beta}$ with $\beta>n$ or $\rho(x)=\exp({-(1+\vert x\vert^2)^{a/2}})$ with  $a>0$. Then there exists a constant $C$ such that for all $t>0$ and $x,y\in \R^n$ the kernel density $p_t$ satisfies  
%$V(x)=(1+\vert x\vert ^2) ^{\beta/2}, x\in \R^n$. Then $V$ satisfies the condition of the lemma
%(\ref{Lemtechn}) for any $\beta >n$ for any $p\in (1,2-\frac{n}{\beta})$.  Hence $V$  is a Lyapunov function for $L$ in 
%${\mathcal L}^p(\mu)$ for such $p$.
%The Lyapunov constant ${\kappa}_{\beta}$ satisfies  $c_{\beta}\leq \beta n+\vert \beta-2\vert$ and for any $t>0, x, y\in \R^n$:
$$
p_{t}(x,y)\leq  {{\frac{C}{t^{n/2}}}\,{ e^{Ct}}\rho^{-1/2}(x)}\,\rho^{-1/2}(y).
$$

%%%%
\ecor

In the next section we shall improve the bound on the kernel density in the case of the measure with density $\rho(x)=\exp({-(1+\vert x\vert^2)^{a/2}})$ with  $a>1$ ; for that purpose  we shall use a Lyapunov function $V$ which will be now in $\mathcal L^2(\mu)$.

%%%%%%%%%%%%%%%%%%%%%%%%%%%%%%%%%%%%%% 
%%%%%%%%%%%%%%%%%%%%%%%%%%%%%%%%%%%%%%
 \section{The measures on $\R$ between exponential and Gaussian}%++%
\label{sec-app}
  In this section we shall prove that the  weighted Nash inequality~\eqref{ineg:NashPhi} holds  with power functions $\phi$ and $\mathcal L^2$ weights $V$ for the semigroups on $\R$ with the invariant measure  $\exp(-|x|^a) dx$. Again for convenience we will study $\exp(-(1+x^2)^{a/2}) dx$ instead of $\exp(-|x|^a) dx$.

  The analysis made here would make no difference if one would work on $\R^n$, except for the values  of the involved constants.  We shall let 
  $$
  T(x)= (1+x^2)^{1/2},
  $$
  and for $a>0$   the probability measure 
  $$
  d\mu_a= C_a e^{-T^a} dx,
  $$ 
  where $C_a$ is the normalizing constant. 
  
  We are dealing with the Sturm-Liouville operator 
  $$
  Lf= f'' -a T^{a-1} T' f',
  $$
  which is symmetric  (and even self adjoint) with respect to the probability measure $\mu_a$. We let  $\rho_a$ denote the density function of the measure $\mu_a$ with respect to the Lebesgue measure, that is
  $$
  \rho_a= \exp(-T^a).
  $$

   In this case, $f$ is in the domain of the Dirichlet form as soon as $f'\in \cL^2(\mu_a)$ and 
  $$
  \cE_{\mu_a}(f,f)= \int f'^2 d\mu_a=-\int fLfd\mu_a.
  $$
  
  \medskip
  
  We shall not pay too much attention to the values of the constants which may be far from being optimal.

  \blem
  \label{lem-liapu}
  For all  $a >0$ and $\beta \in \dR$  the function
  \begin{equation}
  \label{def-weight}
  V= \rho_a^{-1/2} T^{-\beta}=\exp\PAR{\frac{T^a}{2}}T^{-\beta}
  \end{equation}
is a Lyapunov function; moreover  $V\in \cL^2(\mu_a)$ as soon as $\beta>1/2$. 
\elem

 \bpf
First observe that $V$ is positive, and is a Lyapunov function with constant $c$ if and only if
$$
L(\log V)+ (\log V)'^2 \leq c.
$$
But, with $T = T(x)$,
\begin{eqnarray*}
  L(\log V) +(\log V)'^2
&  \!\!=\!\! &
  \frac{a}{2} (a-1) T^{a-2} T'^2 - \frac{a^2}{4} T^{2a-2} T'^2 + \beta (\beta+1) \frac{T'^2}{T^2} + \frac{a}{2} T^{a-1} T''  - \beta \frac{T''}{T} \\
& \!\!=\!\! &
\frac{a}{4}  T^{a-4}  \big( 2(a-1)x^2 - a T^a x^2 + 2 \big) + \beta ( \beta+1) x^2 T^{-4} - \beta T^{-4}
 \end{eqnarray*}
since $T'(x) = x \, T(x)^{-1}$ and $T'' = T(x)^{-3}.$ Now for all $a >0$ the bracket tends to $0$ as $\vert x \vert$ tends to $+\infty$ and for all $\beta$ the last two terms go to $0$, so the continuous map $L(\log V)+ (\log V)'^2$ is bounded from above on $\dR$.
\epf
  
  \bigskip
  
   The  first basic result  is the following
 
  \blem\label{lem:f0=0}  For all $a \geq 1$ and $\beta>0$ there exists a constant $C= C(a, \beta)$ such that,
  for all  smooth and compactly supported  functions $f$ such that $f(0)=0$,
  \begin{enumerate}[(i)]
  \item 
 $$
   \int f^2 d\mu_a \leq C \cE_{\mu_a}(f,f),
  $$
  \item 
 $$
  \int f^2 d\mu_a \leq C \cE_{\mu_a}\PAR{f,f}^\gamma\PAR{\int|f|V d\mu_a}^{2(1-\gamma)}
$$
  where $V$ is the weight given by~\eqref{def-weight} and
  $\displaystyle  \gamma = 1-2\frac{a-1}{3(a-1)+2\beta} \in \big( \frac{1}{3}, 1 \big].$
\end{enumerate}  
   \elem
   
    \bpf 
    We shall let $C$ denote diverse constants depending only on $a$ in the proof of $(i)$, and only on $a$ and $\beta$ in the proof of $(ii)$.
    
    For $x>0$ we let  $q(x)= \int_x^\infty d\mu_a(y)$.  The argument will be based on the following classical estimate (see for instance~\cite[Corollaire 6.4.2]{logsob}):
  \beq
  \label{queue-estim}
  q(x) \leq C \frac{\rho_a(x)}{T(x)^{a-1}}.
  \eeq  
To prove $(i)$, and for $f$ satisfying $f(0)=0$,  we write
  $$
  \int_0^\infty f^2 d\mu_a = 2\int\int_{t=0}^x  f(t)f' (t) d\mu_a(x) dt=2 \int_0^\infty f(t)f'(t)q(t) dt.
  $$
 But, by~\eqref{queue-estim},  we have the upper bound $q(t) \leq C \rho_a(t)$ since $a, T \geq 1$, so that 
   $$
  \int_0^\infty f^2 d\mu_a \leq C \, \|f\|_2 \, \cE_{\mu_a}(f,f)^{1/2}
  $$ 
by the Cauchy-Schwarz inequality.  A similar result holds for the integral on $(-\infty, 0]$, which gives $(i).$
 
 \medskip
 
 Let us now prove $(ii)$ for $a>1$, since for $a=1$ it amounts to $(i)$. Without loss of generality, we assume that $f$ is non-negative. Then 
   $$
   \int_0^\infty f^2 d\mu_a = \int_0^\infty f^2\I_{\BRA{ \frac{f}{\NRM{f}}_2\leq VZ^{-1/2}}}d\mu_a + \int_0^\infty f^2\I_{\BRA{ \frac{f}{\NRM{f}}_2> VZ^{-1/2}}} d\mu_a
   $$
where $Z$ is a positive constant to be chosen later on.
The first term is bounded from above by  $\NRM{f}_2 Z^{-1/2}\int f V d\mu_a$. Then we write the second one as 
 \begin{equation}\label{2eterme}
 \int_0^\infty f^2 \, \I_{\BRA{ \frac{f}{\NRM{f}_2}> VZ^{-1/2}}} d\mu_a 
 =
  2\int_{0}^\infty f(t) \, f'(t) \, [ \int_{t}^\infty \I_{\BRA{ \frac{f(x)}{\NRM{f}_2}> V(x)Z^{-1/2}}} d\mu_a(x)] \, dt
\end{equation}
  by writing $f^2(x) = 2\int_0^x f(t)f'(t) dt$. We bound the inner integral in the following two ways.
  
On the one hand
\begin{equation}\label{1eint}
   \int_{t}^\infty \I_{\BRA{ \frac{f(x)}{\NRM{f}_2}> V(x)Z^{-1/2}}} d\mu_a(x)\leq  \int_{t}^\infty  d\mu_a(x)=q(t) \leq C \rho_a(t)  T(t)^{1-a} 
\end{equation}
 according to~\eqref{queue-estim}.
 
 On the other hand the map $y \mapsto e^{y/2} y^{-\beta}$ is decreasing on $(0, 2 \beta]$ and then increasing, and $T \geq 1$; hence $V$ is increasing on $(0, +\infty)$ if $2 \beta \leq 1$, and it is decreasing on $(0, \sqrt{4 \beta^2 -1}]$ and then increasing otherwise. Hence, in any case, 
  there exists $C$ such that $V(x)\geq C V(t)$ for all $x\geq t > 0$. Hence 
\begin{equation}\label{2eint}
   \int_{t}^\infty \I_{\BRA{ \frac{f(x)}{\NRM{f}_2}> V(x)Z^{-1/2}}} d\mu_a(x)\leq  \int_{t}^\infty \I_{\BRA{ \frac{f(x)}{\NRM{f}_2}> C V(t)Z^{-1/2}}} d\mu_a(x)\leq \frac{Z}{C^2 V^2(t)} = \frac{Z}{C^2} \rho_a(t) T(t)^{2 \beta}
\end{equation}
by the Markov inequality. 

Therefore 
$$
  \int_{t}^\infty \I_{\BRA{ \frac{f(x)}{\NRM{f}_2}> V(x)Z^{-1/2}}} d\mu_a(x) \leq C\rho_a (t) \, \min\BRA{T(t)^{1-a},T^{2\beta}(t)Z}.
$$

Now, since $a+2\beta-1>0$ and $T$ is increasing, then for any $Z\in(0,1]$ there exists $t_0$ such that 
$T(t_0)^{a+2\beta-1}=1/Z$, that is, $T(t_0)^{1-a}=T(t_0)^{2\beta}Z$. We split the integral in~\eqref{2eterme} into two parts, according to $t\geq t_0$ or not, and obtain 
\begin{eqnarray*}
\int_0^\infty f^2 \, \I_{\BRA{ \frac{f}{\NRM{f}_2}> VZ^{-1/2}}} d\mu_a
&\leq&
 C Z \int_0^{t_0} \vert ff'  \vert \, T^{2\beta} \, d\mu_a+C\int_{t_0}^\infty  \vert  ff'  \vert \, T^{1-a} \, d\mu_a \\
 & \leq & 
 C Z T^{2 \beta} (t_0) \int_0^{t_0} \vert ff'  \vert \, d\mu_a+ C T^{1-a} (t_0) \int_{t_0}^\infty  \vert ff'  \vert \, d\mu_a 
\end{eqnarray*}
since $\beta>0$ and $1-a<0$. Moreover $Z T^{2\beta}(t_0) = T(t_0)^{1-a}$, so 
$$
\int_0^\infty f^2\I_{\BRA{ \frac{f}{\NRM{f}_2}> VZ^{-1/2}}} d\mu_a\leq CT(t_0)^{1-a}\NRM{f}_2 \cE_{\mu_a}(f,f)^{1/2}.
$$
by the Cauchy-Schwarz inequality. 

\medskip

In the end we have obtained the bound
$$
\|f\|_2\leq C\SBRA{Z^{-1/2}\int fVd\mu_a+Z^{\frac{1-a}{1-a-2\beta}} \cE_{\mu_a}(f,f)^{1/2}}
$$
for all $0 < Z \leq 1.$

\smallskip

If ${\int fVd\mu_a}\leq{ \cE_{\mu_a}(f,f)^{1/2}}$ then we  choose 
$$
Z=\PAR{\frac{\int fVd\mu_a}{ \cE_{\mu_a}(f,f)^{1/2}}}^{\frac{2(1-a-2\beta)}{3(1-a)-2\beta}} \in (0,1]
$$
to get the inequality
$$
  \int_0^\infty f^2 d\mu_a \leq C \cE_{\mu_a}\PAR{f,f}^\gamma\PAR{\int f V d\mu_a}^{2(1-\gamma)},
$$
where $\gamma = \PAR{a-1+2\beta}/\PAR{3(a-1)+2\beta}$. The same estimate holds on $(-\infty,0]$ which gives $(ii)$. 

If now ${ \cE_{\mu_a}(f,f)^{1/2}}\leq{\int fVd\mu_a}$, then, by $(i)$,
$$
\int f^2 d\mu_a \leq C \cE_{\mu_a}(f,f)= C\cE_{\mu_a}(f,f)^\gamma \cE_{\mu_a}(f,f)^{1-\gamma} \leq C\cE_{\mu_a}(f,f)^\gamma \, \Big( \int fVd\mu_a \Big)^{2(1 - \gamma)}
$$
for all $0 \leq \gamma \leq 1$, which gives $(ii)$.
\epf
  
  \brmq The first point of Lemma \ref{lem:f0=0} is only based on the tail estimate $q(t) \leq C \rho_a (t)$, so holds for all measures $d\mu = \rho \, dx$ such that  $q(x)\leq C \rho(x)$ where $q(x) = \mu([x, +\infty)).$ In particular such probability measures $\mu$ satisfy a spectral gap inequality 
  $$
  \|f\|_2^2 \leq \PAR{\int f  d\mu}^2 + C \cE_{\mu}(f,f)
  $$
by applying $(i)$ to $f - f(0)$,   since 
  $$
{\rm{Var}}_{\mu}(f):=  \int f^2 d\mu -\PAR{\int f d\mu}^2 \leq \int (f-c)^2 d\mu
  $$
  for all constants $c$, and in particular for $c = f(0).$

In fact the probability measure $\mu_a$ is log-concave on $\dR$ and, according to the Bobkov Theorem (see~\cite{bobkov99}), all log-concave measures on $\dR^n$ satisfy a Poincar\'e inequality. Note that a proof of this result is given in~\cite{guillin-al-08-1} by using the Lyapunov function $W=e^{\gamma T^a}$ for a $\gamma>0$.   
  \ermq

  \brmq The condition $a \geq 1$ is crucial in this proof of Lemma~\ref{lem:f0=0}. The  second point is obtained for all $\beta >0$. For $\beta \leq 0$ we may use the bound~\eqref{2eint} with $T(t)^{2 \beta} \leq 1$, but not~\eqref{1eint}; then we choose $Z = (\int f V d\mu_a \cE_{\mu_a}(f,f)^{-1/2} )^{2/3}$ to obtain $(ii)$ with $\gamma= 1/3$. Observe that the best bound is obtained for $\beta=0$, for which we have the following general bound.
  \ermq
  
  \brmq
 Let $\mu$ be a probability measure on $\R$, with  a density $\rho(x)$ increasing on $(-\infty,0)$ and decreasing on $(0,\infty)$ and let $V= \rho^{-1/2}$. Then 
  $$
  \|f\|_2 \leq  \Big( \frac{27}{2} \Big)^{1/3} \PAR{\int|f|V d\mu}^{1/3}\cE_{\mu} (f,f)^{1/3}
  $$
for all smooth functions such that $f(0)=0$. 
 The proof follows the argument of Lemma~\ref{lem:f0=0}, by using the bound~\eqref{2eint} but not~\eqref{1eint}.   It gives a Nash inequality with $\phi(x) = 2 x^{3/2} /27$ on $(0, +\infty)$, so that $1/\phi$ is integrable at infinity. However, besides the restriction $f(0)=0$ which will be removed below only for $a> 3$ (with $\beta =0$),  it does not give any upper bound on the density, as in Corollary~\ref{coro-sect2}, since $V$ is not in $\cL^2(\mu)$. 
  \ermq
  
  The restriction $f(0)=0$ is removed by the following
  
  \blem
  \label{lem:comparaison} 
  Given the measure $d\mu_a=C_a \exp(-T^a) dx$ with $a>0$ and the weight function 
  $$
  V= \exp(T^a/2)T^{-\beta}
  $$ with 
  $$
  \beta> \frac{3-a}{2},
  $$
   then there exist  $\theta\in(0,1)$ and constant $C$ such that 
    $$
  \int |f-f(0)|V d\mu_a \leq C \SBRA{\int |f|V d\mu_a + \PAR{\int |f|V d\mu_a}^{1-\theta}\cE_{\mu_a}(f,f)^{\theta/2}}
  $$ 
for all nonnegative smooth compactly supported $f$ on $\R$.
  \elem
  
  \brmq\label{rem: comparaison}
  For $\beta > 3/2$ then all $\theta \in (2/3, 1)$ are admissible.
  \ermq
  
  \bpf  
  In the proof we shall let $C$ denote diverse constants which depend only on $a, \beta$ and a parameter $\alpha$ to be introduced later on. We start by writing
  \begin{equation}
  \label{eq-f1}
  \int\ABS{f-f(0)} V d\mu_a\leq \int \ABS{f}Vd\mu_a+\ABS{f(0)}\int Vd\mu_a.
  \end{equation}
For convenience  we let
   $$
   U= \int|f|V d\mu.
   $$
For any $\alpha>0$, and any  $x\in\dR$,  write
  $$
  |f^\alpha(x)-f^\alpha(0)|= \alpha \, |\int_0^x f^{\alpha-1} f' dx|\leq C \ABS{\int_0^x |fV|^{\alpha-1} |f'| \frac{1}{\rho_a V^{\alpha-1}} d\mu_a}.
  $$
By the H\"older inequality, for  any   $p,q,r>1$ such that ${1}/{p}+{1}/{q}+ {1}/{r}=1$, then
  $$
  \ABS{\int \!fgh d\mu_a}\leq \|f\|_p\|g\|_q\|h\|_r.
  $$
For $q=2$, $p=1/(\alpha-1)$ and $r=2/(3-2\alpha)$ with $\alpha\in(1,3/2)$ this gives
  $$
  |f^\alpha(x)-f^\alpha(0)|\leq C \,{U}^{\alpha-1}\cE_{\mu_a}(f,f)^{1/2} K^\alpha(x),
  $$
 where 
  $$
  K(x)= \ABS{\int_0^x \frac{\rho_a(t) ^{1-r}}{V(t)^{r(\alpha-1)}} dt}^{\frac{1}{r\alpha}}.
  $$
 
  Then
  $$
  \ABS{f(0)}\leq\ABS{f(x)}+\ABS{f^\alpha(0)-f^\alpha(x)}^{1/\alpha}
  $$
for all $x$ since $\alpha  \geq 1$,  so
  $$
  \ABS{f(0)} \leq C\SBRA{|f(x)|+ U^{1-1/\alpha}\cE_{\mu_a}(f,f)^{1/(2\alpha)} K(x)},
  $$
and then
  \begin{equation}
  \label{eq-f0}
  |f(0)|\int V d\mu_a \leq C\SBRA{U+ U^{1-1/\alpha}\cE_{\mu_a}(f,f)^{1/(2\alpha)}\int K V d\mu_a}.
  \end{equation}
Let us prove that   $\int K V d\mu_a$ is finite. By the definition~\eqref{def-weight} of $V$ and~\cite[Corollaire 6.4.2]{logsob} for instance, one has 
  \begin{equation}
  \label{eq-majq}
  K(x)= \ABS{\int_0^x e^{\frac{T^a}{2}\PAR{\frac{3}{2}r-1}}T^{\beta r( \alpha-1)}dt}^{\frac{1}{r\alpha}}
  \leq
 C\exp\PAR{\frac{T^a(x)}{2}}T^{d}(x),
  \end{equation}
  with 
  $$
  d= \beta\PAR{1-\frac{1}{\alpha}}-\frac{a-1}{r \alpha}.
  $$ 
  In fact the two quantities in~\eqref{eq-majq} are equivalent when $\ABS{x}$ is large.
  
Hence  $ K V\rho_a\leq CT^{d-\beta}$, so the  integral $\int K V d\mu_a$ is convergent as soon as $d-\beta<-1$, that is,
   $$
  \alpha < 1 + \frac{1}{a} \Big( \beta - \frac{3-a}{2} \Big).
  $$
  Hence, if $\beta > (3-a)/2$, then any $\displaystyle 1 < \alpha < \min \Big\{ \frac{3}{2}, 1 + \frac{1}{a} \big( \beta - \frac{3-a}{2} \big) \Big\}$ satisfies all conditions, so that $\int K V d\mu_a<\infty$.  Then
  $$
  \int|f-f(0)| Vd\mu_a \leq C\SBRA{U+ U^{1-1/\alpha}\cE_{\mu_a}(f,f)^{1/(2\alpha)}}
  $$
 by~\eqref{eq-f1} and~\eqref{eq-f0}. This proves Lemma~\ref{lem:comparaison}  with $\theta= 1/\alpha$. 
 \epf
  
  \brmq
 The argument is only based on the fact that the function 
  $$
  K(x)= \ABS{\int_0^x \frac{\rho_a^{1-r}}{V^{r (\alpha-1)}} dt}^{\frac{1}{r\alpha}}
  $$ 
  satisfies 
  $$
  \int K \rho_a V dx <\infty.
  $$
  In particular, in the limiting case when $\beta =0$ and $V(x)= \rho_a^{-1/2}$, this amounts to
  $$
  \int_0^\infty [\int_0^x \rho_a(t)^{-\alpha/(3-2\alpha)} dt]^{(3-2\alpha)/2\alpha} \rho_a^{1/2} (x)dx < \infty,
  $$
  that is, $a> \displaystyle \frac{3}{3 - 2 \alpha}$ (see again~\cite[Corollaire 6.4.2]{logsob} for instance). In turn this holds for an $\alpha \in (1, 3/2)$ if and only if $a>3.$
   \ermq
  
  \brmq
The two fundamental lemmas are based on the two estimates~\eqref{queue-estim} and~\eqref{eq-majq}. These are basic estimates when proving that a probability measure on $\dR$ satisfies a Poincar\'e or a logarithmic Sobolev inequalities, as explained in~\cite[Section 6.4]{logsob}.   
  \ermq
  
  Collecting lemmas \ref{lem:f0=0} and \ref{lem:comparaison}, we get the following main result:
  \bthm \label{thm:main}  On $\R$,  let us consider the measure
   $$
   d\mu_a(x)= C_a\exp(-T^a) dx
  $$
   with $T(x)= (1+|x|^2)^{1/2}$,  and the weight function 
    $$
    V= \exp\PAR{\frac{T^a}{2}} T^{-\beta}
    $$
     with $a >1$ and  $\beta \in \dR$.
 Then there exist  $C$ and $\lambda   \in (0, 1)$ such that
  \begin{equation}
  \label{eq-derthm}
  \|f\|_2^2 \leq C\SBRA{\PAR{\int |f|V d\mu_a}^2 + \PAR{\int |f|V d\mu_a}^{2(1-\lambda)}\cE_{\mu_a}(f,f)^{\lambda}}
   \end{equation}
   for all functions $f$ . 
  \ethm
 
  \bpf
The space of smooth functions with compact support is dense in the domain of $L$, so it is enough to consider the case when $f$ is smooth and compactly supported. Also, without loss of generality, we may assume that $f$ is nonnegative. Here again $C$ will denote diverse constants depending on the parameters $a$ and $\beta$ and a parameter $\theta$ to be introduced later on.
   
  One has, 
  $$
  \|f\|_2^2 \leq \PAR{\int fd\mu_a}^2 + \int |f-f(0)|^2 d\mu_a.
  $$
 The weight $V$ is bounded from below by a positive constant, so
  \begin{equation}
  \label{eq-inthm}
  \|f\|_2^2 \leq C \PAR{\int fVd\mu_a}^2 + \int |f-f(0)|^2 d\mu_a.
  \end{equation}
  Let now $U= \int f V d\mu_a$ and $U_0= \int |f-f(0)| V d\mu_a$ and assume $\beta >0$.  By Lemma \ref{lem:f0=0},  applied to the function  $f- f(0)$, one has
  \begin{equation}
  \label{eq-inthm2}
  \int |f-f(0)|^2 d\mu_a \leq C{\cE_{\mu_a}(f,f)}^{\gamma} U_0^{2(1-\gamma)}, 
  \end{equation}
   where 
  $$
  \gamma = 1- 2\frac{a-1}{3(a-1)+ 2 \beta}.
  $$
But, if moreover $\beta > (3-a)/2$, by Lemma \ref{lem:comparaison} there exists $\theta\in (0,1)$ such that
  $$
  U_0 \leq C \SBRA{U+ U^{1-\theta}\cE_{\mu_a}(f,f)^{\theta/2} },
  $$
so that
  $$
  \int |f-f(0)|^2 d\mu_a \leq C {\cE_{\mu_a}(f,f)}^{\gamma} \SBRA{U^{2(1-\gamma)}+U^{2(1-\theta)(1-\gamma)}  \cE_{\mu_a}(f,f)^{\theta(1-\gamma)}}
  $$
by~\eqref{eq-inthm2}. Hence, by~\eqref{eq-inthm},
  $$
  \|f\|_2^2 \leq CU^2\SBRA{1+\PAR{\frac{{\cE_{\mu_a}(f,f)}}{U^2}}^{\gamma}+\PAR{\frac{{\cE_{\mu_a}(f,f)}}{U^2}}^{\gamma+\theta(1-\gamma)}}
  \leq C\SBRA{U^2+{{{\cE_{\mu_a}(f,f)}}^{\lambda}{U}}^{2(1-\lambda)}}
  $$
 if $\lambda = \gamma + \theta (1 - \gamma) \in (0,1)$. This concludes the argument  for $\beta > \max(0, \frac{3-a}{2}).$
 
 \smallskip
 
 Then, since $V$ is decreasing in $\beta$, then~\eqref{eq-derthm} holds for all real $\beta.$
     \epf
  
  \brmq 
 
We are restricted to $a>1$, since for $a=1$ then only $\lambda=1$  is admissible; this gives a useless inequality for our purpose, which is even weaker than the Poincar\'e inequality. 

  According to Lemma~\ref{lem:comparaison} and Remark~\ref{rem: comparaison} the larger $\beta$ is, the smaller the weight $V$ is, and the larger exponent $\lambda$ of the Dirichlet form has to be in~\eqref{eq-derthm}; on the contrary, the smaller $\beta$ is $(> 3/2)$, the smaller exponent $\lambda$  we can take.
  \ermq 
 
 \bigskip
 
 We can now illustrate the abstract method  of section~\ref{sec-prems} by obtaining the following pointwise bounds on the Markov semigroup associated to $L$, which bring new information on this semigroup for small time:
  
 \bcor\label{cor-main}
 Let $a >1$ and let $(P_t)_{t\geq0}$ be the Markov generator on $\R$ with generator
  $$
  Lf= f'' - aT^{a-1}T' f',
  $$ 
  and reversible measure $d\mu_a(x) = \rho_a(x)dx= C_a\exp(-(1+|x|^2)^{a/2})dx$. 
  
  Then for all real $\beta$ there exists $\delta>0$ 
   and a constant $C$ such that, for all $t$, $P_t$ has a density $p_t$ with respect to the measure $\mu_a$, which satisfies
  $$
  p_t(x,y)\leq  \frac{C e^{Ct}}{t^{\delta }}\,\frac{\rho_a^{-1/2}(x)\rho_a^{-1/2}(y)}{(1+|x|^2)^{\beta/2}(1+|y|^2)^{\beta/2}}
  $$
for almost every $x,y\in\dR$.

  Moreover, the spectrum of $-L$ is discrete and its eigenvalues $(\lambda_n)_{n\in\N}$ satisfy the inequality 
  $$
  \sum_n e^{-\lambda_nt} \leq \frac{C e^{Ct}}{ t^{\delta}}
  $$
for all $t>0$.
  \ecor

  \bpf
Letting $C$ and $\lambda \in (0, 1)$ be defined as in Theorem~\ref{thm:main}, by the inequality~\eqref{eq-derthm} the Dirichlet form $\cE_{\mu_a}$ satisfies a weighted Nash inequality with weight $V = \exp (T^a /2) T^{-\beta}$  and rate function 
  $$
  \phi(x) = C^{-1/\lambda} (x-C)^{1/\lambda}
  $$
on $(C, +\infty)$.  Moreover the weight $V$ is a Lyapunov function with constant $c>0$ by Lemma~\ref{lem-liapu}, it is in $\cL^2(\mu_a)$ if $\beta > 1/2$ and hypothesis~\eqref{eq-hypom} of Theorem~\ref{thm:BornesPt} holds since $\lambda <1$. Hence, by Corollary~\ref{coro-sect2} and for diverse constants $C= C(a, \beta, \lambda)$, for all $t>0$ the operator $P_{2t}$ has a density $p_{2t}$ with respect to $\mu_a$, which satisfies
 $$
  p_{2t}(x,y)
  \leq  
  C(1+ t^{\frac{-2 \lambda}{1-\lambda}} ) e^{2 ct}\frac{\rho_a^{-1/2}(x)\rho_a^{-1/2}(y)}{(1+|x|^2)^{\beta/2}(1+|y|^2)^{\beta/2}}
\leq  Ct^{\frac{-2 \lambda}{1-\lambda}} e^{2ct} \frac{\rho_a^{-1/2}(x)\rho_a^{-1/2}(y)}{(1+|x|^2)^{\beta/2}(1+|y|^2)^{\beta/2}}.
  $$
  This proves the first statement  for $\beta >1/2$, with $\delta=2 \lambda/(1-\lambda) >0$, and then for any $\beta$.  
  
 The second statement on the trace of $P_t$ is obtained by letting any $\beta > 1/2$ in the upper bound on $p_t (x,x)$ and integrating. 
  \epf
 \medskip
 
For $\beta>1/2$,  the non-uniform bound implies that $P_t$ is Hilbert-Schmidt but we do not recover the Orlicz hypercontractivity result of \cite{bcr1}. This is not surprising since  in fact no bound such as $K(t)V(x)V(y)$ can imply hypercontractivity of more generally Orlicz hypercontractivity.

 \brmq
  \label{rm-last} The same method, with $V=1$, leads to a (non weighted) Nash inequality for $\mu_a$ with $a >1$, with rate function
   $$
   \phi(x)= C \, x \, (\log x)^{2(1-1/a)}
   $$
    on an interval $(M, \infty)$. By Theorem~\ref{thm:BornesPt} this implies  that the semigroup is ultracontractive as soon as $1/\phi$ is integrable at infinity, that is, for  $a>2$, hence recovering a partial result of~\cite{kkr}.
      \ermq

  \brmq 
    Observe in Corollary~\ref{cor-main} that there is no optimal $\beta$, that is, no optimal bound on $p_t(x,x)$ of the form $C(t)\rho(x)^{-1/2}T(x)^{-\beta}$.
  So one could look for an optimal bound on $p_t(x,x)$  such as $C(t) \, \rho^{-\lambda}(x)$ for a $\lambda\in(0,1/2)$. It is not the case in the Gaussian case when $a=2$: in this case the optimal bound is $C(t) \, \exp \big( \vert x \vert^2 /(1+ e^{2t}) \big)$, hence of the form $C(t) \, \rho(x)^{- \lambda (t)}$ with  $\lambda (t) < 1/2$; it is even an equality, see~\eqref{boundNt}.
  
  Also for $1<a < 2$ it seems that $p_t(x,x)$ can not be bounded by $C(t) \rho^{-\lambda}(x)$ for $\lambda < 1/2$. Indeed, for the weight $V= \exp(\lambda T^a)$ with $\lambda<1/2$, our method leads to a weighted Nash inequality with rate function 
$$
\phi(x)= C(a,\lambda) \, x \, (\log x)^{2(1-1/a)}
$$ 
 on an interval $(M, \infty)$, where $\lambda$ appears only in the value of the constant $C(a,\lambda)$. Apart from the values of the constants, this is not better than the inequality obtained in Remark~\ref{rm-last} with $V=1$, and again this is not enough to obtain any bound on the density $p_t(x,y)$ by lack of integrability of $1/\phi$. Now we do not know whether a bound such as  $C(t) \, \rho(x)^{- \lambda (t)}$ with  $\lambda (t) < 1/2$ could be optimal, but we strongly doubt about it.

Again from this point of view the Gaussian case appears as a particular case, being a critical case as regards the two points of view of ultracontractivity and non-uniform bounds; in this case, and in this case only, one may do better, and Gaussian Nash inequalities are under study in a work in progress.
\ermq

\noindent
{\bf Acknowledgements.} We would like to thank F.-Y. Wang for pointing out that Theorem~\ref{thm:BornesPt} is strongly related to~\cite[Theorem~3.3]{wang02}.

This research was supported in part by the ANR project EVOL. The third author thanks the members of UMPA at the Ecole Normale Sup\'erieure de Lyon for their kind hospitality.

\bigskip
\noindent

%\bibliographystyle{plain}
%\bibliography{WNI-biblio}

\begin{thebibliography}{10}

\bibitem{logsob}
C.~An{\'e}, S.~Blach{\`e}re, D.~Chafa{\"{\i}}, P.~Foug{\`e}res, I.~Gentil,
  F.~Malrieu, C.~Roberto, and G.~Scheffer.
\newblock {\em Sur les in\'egalit\'es de {S}obolev logarithmiques}, volume~10
  of {\em Panoramas et Synth\`eses}.
\newblock Soci\'et\'e Math\'ematique de France, Paris, 2000.

\bibitem{bakrystflour}
D.~Bakry.
\newblock L'hypercontractivit\'e et son utilisation en th\'eorie des
  semigroupes.
\newblock In {\em Lectures on probability theory ({S}aint-{F}lour, 1992)},
  Lecture Notes in Math. 1581, pages 1--114. Springer, Berlin, 1994.

\bibitem{guillin-al-08-1}
D.~Bakry, F.~Barthe, P.~Cattiaux, and A.~Guillin.
\newblock A simple proof of the {P}oincar\'e inequality for a large class of
  probability measures including the log-concave case.
\newblock {\em Electron. Commun. Probab.}, 13:60--66, 2008.

\bibitem{bbg1}
D.~Bakry, F.~Bolley, and I.~Gentil.
\newblock Dimension dependent hypercontractivity of gaussian kernels.
\newblock {\em Preprint}, 2010.

\bibitem{guillin-al-08-2}
D.~Bakry, P.~Cattiaux, and A.~Guillin.
\newblock Rate of convergence for ergodic continuous {M}arkov processes:
  {L}yapunov versus {P}oincar\'e.
\newblock {\em J. Funct. Anal.}, 254(3):727--759, 2008.

\bibitem{BCLS}
D.~Bakry, T.~Coulhon, M.~Ledoux, and L.~Saloff-Coste.
\newblock Sobolev inequalities in disguise.
\newblock {\em Indiana Univ. Math. J.}, 44(4):1033--1074, 1995.

\bibitem{bcr1}
F.~Barthe, P.~Cattiaux, and C.~Roberto.
\newblock Interpolated inequalities between exponential and {G}aussian,
  {O}rlicz hypercontractivity and isoperimetry.
\newblock {\em Rev. Mat. Iberoam.}, 22(3):993--1067, 2006.

\bibitem{bcr07}
F.~Barthe, P.~Cattiaux, and C.~Roberto.
\newblock Isoperimetry between exponential and {G}aussian.
\newblock {\em Electron. J. Probab.}, 12:no. 44, 1212--1237, 2007.

\bibitem{bcs}
A.~Bendikov, T.~Coulhon, and L.~Saloff-Coste.
\newblock Ultracontractivity and embedding into {$L^\infty$}.
\newblock {\em Math. Ann.}, 337(4):817--853, 2007.

\bibitem{bobkov99}
S.~G. Bobkov.
\newblock Isoperimetric and analytic inequalities for log-concave probability
  measures.
\newblock {\em Ann. Probab.}, 27(4):1903--1921, 1999.

\bibitem{bobkov-gotze}
S.~G. Bobkov and F.~G{\"o}tze.
\newblock Exponential integrability and transportation cost related to
  logarithmic {S}obolev inequalities.
\newblock {\em J. Funct. Anal.}, 163(1):1--28, 1999.

\bibitem{cks}
E.~A. Carlen, S.~Kusuoka, and D.~W. Stroock.
\newblock Upper bounds for symmetric {M}arkov transition functions.
\newblock {\em Ann. Inst. H. Poincar\'e Probab. Statist.}, 23(2,
  suppl.):245--287, 1987.

\bibitem{carlen-loss}
E.~A. Carlen and M.~Loss.
\newblock Sharp constant in {N}ash's inequality.
\newblock {\em Internat. Math. Res. Notices}, (7):213--215, 1993.

\bibitem{coulhon}
T.~Coulhon.
\newblock Ultracontractivity and {N}ash type inequalities.
\newblock {\em J. Funct. Anal.}, 141(2):510--539, 1996.

\bibitem{davies}
E.~B. Davies.
\newblock {\em Heat kernels and spectral theory}, volume~92 of {\em Cambridge
  Tracts in Mathematics}.
\newblock Cambridge University Press, Cambridge, 1990.

\bibitem{ggm06}
I.~Gentil, A.~Guillin, and L.~Miclo.
\newblock Modified logarithmic {S}obolev inequalities in null curvature.
\newblock {\em Rev. Mat. Iberoam.}, 23(1):235--258, 2007.


\bibitem{gross75}
L.~Gross.
\newblock Logarithmic {S}obolev inequalities.
\newblock {\em Amer. J. Math.}, 97(4):1061--1083, 1975.

\bibitem{kkr}
O.~Kavian, G.~Kerkyacharian, and B.~Roynette.
\newblock Quelques remarques sur l'ultracontractivit\'e.
\newblock {\em J. Funct. Anal.}, 111(1):155--196, 1993.

\bibitem{kolmogorov}
A.~Kolmogorov, S.~Fomine, and V.~M. Tihomirov.
\newblock {\em El\'ements de la th\'eorie des fonctions et de l'analyse
  fonctionnelle}.
\newblock \'Editions Mir, Moscow, 1974.

\bibitem{latala}
R.~Lata{\l}a and K.~Oleszkiewicz.
\newblock Between {S}obolev and {P}oincar\'e.
\newblock In {\em Geometric aspects of functional analysis}, volume 1745 of
  {\em Lecture Notes in Math.}, pages 147--168. Springer, Berlin, 2000.

\bibitem{ledoux2000}
M.~Ledoux.
\newblock The geometry of {M}arkov diffusion generators.
\newblock {\em Ann. Fac. Sci. Toulouse Math. (6)}, 9(2):305--366, 2000.

\bibitem{maheux}
P.~Maheux.
\newblock Nash-type inequalities and decay of semigroups of operators.
\newblock {\em Preprint}, 2010.

\bibitem{muckenhoupt}
B.~Muckenhoupt.
\newblock Hardy's inequality with weights.
\newblock {\em Studia Math.}, 44:31--38, 1972.

\bibitem{nash}
J.~Nash.
\newblock Continuity of solutions of parabolic and elliptic equations.
\newblock {\em Amer. J. Math.}, 80:931--954, 1958.

\bibitem{saloff11}
L.~Saloff-Coste.
\newblock Sobolev inequalities in familiar and unfamiliar settings.
\newblock In {\em Sobolev spaces in mathematics. {I}}, volume~8 of {\em Int.
  Math. Ser. (N. Y.)}, pages 299--343. Springer, New York, 2009.

\bibitem{va-sc-co}
N.~Th. Varopoulos, L.~Saloff-Coste, and T.~Coulhon.
\newblock {\em Analysis and geometry on groups}, volume 100 of {\em Cambridge
  Tracts in Mathematics}.
\newblock Cambridge University Press, Cambridge, 1992.

\bibitem{wang2000}
F.-Y. Wang.
\newblock Functional inequalities for empty essential spectrum.
\newblock {\em J. Funct. Anal.}, 170(1):219--245, 2000.


\bibitem{wang02}
F.-Y. Wang.
\newblock Functional inequalities and spectrum estimates: the infinite measure
  case.
\newblock {\em J. Funct. Anal.}, 194(2):288--310, 2002.


\bibitem{wangbook}
F.-Y. Wang.
\newblock {\em Functional Inequalities, Markov Processes and Spectral Theory}.
\newblock Science Press. Beijing, 2004.

\bibitem{wang-05}
F.-Y. Wang.
\newblock A generalization of {P}oincar\'e and log-{S}obolev inequalities.
\newblock {\em Potential Anal.}, 22(1):1--15, 2005.

\bibitem{yosida}
K.~Yosida.
\newblock {\em Functional analysis}.
\newblock Springer-Verlag, Berlin, 1995.

\end{thebibliography}

\medskip\noindent

\noindent
Institut de Math\'ematiques de Toulouse, UMR CNRS 5219\\
Universit\'e de Toulouse\\
Route de Narbonne\\
31062 Toulouse - France\\
bakry@math.univ-toulouse.fr 

\medskip

\noindent
Ceremade, UMR CNRS 7534\\
Universit\'e Paris-Dauphine\\
Place du Mar\'echal De Lattre De Tassigny\\
75016 Paris - France\\
bolley,gentil@ceremade.dauphine.fr

\medskip

\noindent
Mapmo, UMR CNRS 6628\\
 Universit\'e d'Orl\'eans\\
B\^atiment de math\'ematiques - Route de Chartres\\
45067 Orl\'eans - France\\
patrick.maheux@univ-orleans.fr

  \end{document}